\documentclass[12pt]{article}
\usepackage{graphicx}
\usepackage{amssymb}

\usepackage[dvips]{epsfig}
\usepackage{subfigure}

\begin{document}

\title{Chaos control via feeding switching in intraguild predation systems
}

\author{
Joydev Chattopadhyay$^{\ddagger}$,
Nikhil Pal$^{**}$,\\
Sudip Samanta$^{\ddagger}$,
Ezio Venturino$^{\dagger}$\thanks{This research was partially supported by the project Metodi numerici in teoria delle popolazioni
of the Department of Mathematics ``Giuseppe Peano''.}\\
$^{\ddagger}$ Agricultural and Ecological Research Unit,\\
Indian Statistical Institute\\
203, B. T. Road, Kolkata 700108, India;\\
$^{**}$ Department of Mathematics, Visva-Bharati University,\\
Santiniketan 731235, India.\\
$^{\dagger}$ Dipartimento di Matematica ``Giuseppe Peano'',\\
via Carlo Alberto 10, 10123 Torino,\\
Universit\`{a} di Torino, Italy.\\
Emails: joydev@isical.ac.in, ezio.venturino@unito.it
}
\date{}
\maketitle

\begin{abstract}
In the present paper we study a tri-trophic Lotka-Volterra food web model with omnivory and predator switching.
We observe that if the intensity of predator switching increases the chaotic behavior of the omnivory model will be reduced.
We conclude that the predator switching mechanism enhances the stability of an otherwise chaotic system.
\end{abstract}


\section{Introduction}\label{Introduction}

An interesting and appealing topic in population ecology is the understanding of the diversity and community composition of
the populations which ultimately determine the overall stability of an ecosystem.
The predators have in general a
rigid type of feeding pattern, which may be the choice of a specific prey
or it may be dependent on the abundancy of the prey. In this latter case
looking for resources,
the predator moves toward another prey, which may be present in the same or in another habitat.
This mechanism of preferential predation is named switching, \cite{Tansky1978}.
Many examples of predator feeding switching have been identified in nature,
\cite{Fisher-Piette1934,Murdoch1969}, as well as in laboratory experiments
involving {\it {Notonecta}} and {\it {Ischnura}}, \cite{Lawton1974}.
Several mathematical
models have been proposed with predator switching involving one predator with two
prey species \cite{Holling1961,Teramoto1979,Prajneshu1987,KBJ,KBW,VPB}.

In an omnivory food web, a predator consumes more than one prey species.
In particular, in a tri-trophic food chain with omnivory, top predator predates
both the bottom prey and the intermediate predator \cite{Pimm1978}.
This phenomenon is also termed intraguild predation (IGP), defined
as the feeding of the intermediate predator by a top predator that can also
consume the prey of the intermediate
predator \cite{Polis1989}. Thus, the top predator and the intermediate predator species
are also potentially competitors for the common resource, i.e. the prey at the lowest trophic
level. This fact enhances the chance of switching feeding behavior of the top predator
between the bottom prey and the intermediate predator.

Whether predator switching phenomena act as a stabilizing or destabilizing factors in
a dynamical population system is still debated, see \cite{Kimbrell2004} and
the references therein. A stabilizing effect is found for instance in Ref. \cite{KK,AV}.
Intraguild predation with nonlinear functional responses
had been investigated showing that omnivory could be a stabilizing factor \cite{McCann1997}.
On the other hand, intraguild predation could destabilize a three species Lotka-Volterra
model with linear functional responses, see \cite{Holt1997}.
Furthermore, the Lotka-Volterra model of intraguild predation could
exhibit a limit cycle, coexistence of all species in
IGP system is possible if the intermediate predator is superior at exploitative
competition for the prey, whereas
the top predator gains significantly from its consumption of intermediate predator
\cite{Holt1997b}. Furthermore, the
top predator may become extinct if the predation rate and the conversion efficiency from the
prey to predator are very low, whereas the intermediate predator will disappear
if the predation rate and
the conversion efficiency from prey to predator are high enough. Also recently,
in the very same model, chaotic behavior has been numerically discovered, \cite{Tanabe2005}.
The findings indicate that
the system enters into a chaotic regime from a stable steady state through period doubling
Hopf-bifurcation,
if the consumption rate of top predator on the primary prey increases.

The aim of the present investigation is to incorporate predator
switching in a three species Lotka-Volterra model with linear functional response
and investigate its effect.

The paper is organized as follows. After discussing some preliminary biological background,
we introduce the model in Section 3, and analyse its behavior in the following Section.
The main result, control of the chaotic behavior, is contained in the final Section.

\section{Background}

Predator feeding switching is an important aspect in nature, \cite{Fisher-Piette1934,Holt1997b}.
However, the effect
of predator switching in an omnivory system, also known as intraguild predation,
does not appear to be properly investigated yet. As far as our knowledge goes, no attempt
has been made so far to model
the behavior of such intraguild predation in the presence of predator switching.
In a
tri-trophic food chain or food web with omnivory the top predator consumes more than one prey
population. If the predator population has more than one resource to prey upon,
namely the prey and middle predators,
the predator population predates more heavily either on the preferred or on the most abundant
species, \cite{Murdoch1969}, depending on the specific situation.
Whenever the primary resource declines and while waiting that it recovers,
the predator for its survival feeds on another prey, in the same or another habitat.

Experimental results suggest that predator switching is a very resemble biological phenomena to omnivory
and likely to occur with it \cite{Gismervikl1997}.
However, the predator switching and the omnivory do not yet appear to be studied
simultaneously. In ref. \cite{Holt1997} Holt and Polis have intuitively suggested that
if the adaptive foraging by the omnivore
predator leads to switching between the prey and the intermediate predator,
the system will be stabilized. Our aim is to investigate and substantiate this claim,
by providing a theoretical setting for the effect of predator feeding switching
in omnivory system considering a suitable mathematical model.

\section{Model}

In the three species food chain model proposed by Tanabe and
Namba (2005), \cite{Tanabe2005} chaotic behavior is discovered. We want to consider a modification
of this model, incorporating top predator feeding switching behavior in it.
Let $x$, $y$ and $z$ respectively denote the prey, middle predator and top predator populations sizes.
The resulting system of equations reads as follows
\begin{eqnarray}\label{EQ:eqn1.2}
\begin{array}{lll}
\displaystyle{\frac{d x}{d t}} & = &\displaystyle{ (b_1-a_{11}x)x -a_{12}xy-\frac{a_{13}xz}{1+c\frac{y}{x}}},\\
\displaystyle{\frac{d y}{d t}} & = & \displaystyle{-b_2y+a_{21}xy-\frac{a_{23}yz}{1+c\frac{x}{y}}},\\
\displaystyle{\frac{d z}{d t}} & = &\displaystyle{ -b_3z+\frac{a_{31}xz}{1+c\frac{y}{x}}+\frac{a_{32}yz}{1+c\frac{x}{y}}} ,
\end{array}
\end{eqnarray}
As it occurs for the switching models,
the system (\ref{EQ:eqn1.2}) is not well defined at the origin. To overcome the occurrence of the singularity there,
we simply replace the above equations at $(0,0,0)$ by the following ones:
$$
\frac{dx}{dt}=\frac{dy}{dt}=\frac{dz}{dt}=0.
$$

The meaning of the parameters, all nonnegative, is as follows:
$b_1$ is the prey intrinsic net growth rate, while $a_{11}$ denotes their intraspecific competition for resources.
The parameters $a_{ij}$ with $i<j$ denote the prey consumption rates,
while $a_{ij}$ for $i>j$ represent the corresponding predators' assimilation rates. The natural mortalities of
the intermediate and top predators are respectively $b_2$ and $b_3$.

Therefore, the first equation expresses logistic growth of the prey in absence of predators. When they are present,
the prey are subject to their hunting. In absence of prey, the intermediate predator is bound to starve; in the food
chain it also represents a possible source of food for the top predator. The latter feeds on both prey and intermediate
predators.

The novel feature of this model resides in the last terms of the first two equations and the corresponding last two terms
of the last equation, modelling the feeding switching behavior between the two types of prey for the top predator, following
\cite{Tansky1978}. This mechanism is governed by the relative abundances of the bottom prey and intermediate predators.
The parameter $c$ represents the switching intensity.
This parameter provides also a mean to obtain a continuum of models. In fact,
note that for $c=0$ the system shows no switching behavior, it becomes exactly the Tanabe and Namba three species omnivory system,
\cite{Tanabe2005}, while for $c=1$ the switching function is the same function defined by Tansky,
\cite{Tansky1978}. This parameter thus allows to control the amount of switching feeding that the top predator can use.

\section{The system general behavior}

The steady states that (\ref{EQ:eqn1.2}) can attain are the five points ${E}_k=({x}_k,{y}_k,{z}_k)$, $k=0.\ldots,4$.
The origin $E_0$, representing ecosystem disappearance, fortunately cannot
be stably attained.

The bottom prey-only equilibrium $E_1=(x_1,0,0)$ in which the prey settles at carrying capacity $x_1=b_1a_{11}^{-1}$
is locally asymptotically stable if both predator's mortalities are large enough,
$$
{a_{21}b_1<a_{11}b_2}, \quad \displaystyle{a_{31}b_1<a_{11}b_3}.
$$

The bottom prey can thrive also in presence of either one of the predators,
$$
{E}_2=\left(\frac{b_3}{a_{31}},0,\frac{b_1- a_{11}{x}_2}{a_{13}}\right), \quad
{E}_3=\left(\frac{b_2}{a_{21}},\frac{b_1-a_{11}{x}_3}{a_{12}},0\right).
$$

These equilibria are locally asymptotically stable if the following conditions respectively hold
$$
a_{31}b_2>a_{21}b_3
$$
for $E_2$, i.e. the mortality of the intermediate predator must be higher than the one of the top predator, while for $E_3$ we must have
$$
b_3>\frac{a_{31}{{x_3}^2}}{x_3+cy_3}+\frac{a_{32}{{y_3}^2}}{y_3+cx_3},
$$
expressing the fact that the top predator mortality should be higher than a certain threshold, which is regulated by the other populations
equilibrium values.

\begin{figure}[ht]
\centering
%\begin{tabular}{c}
%\includegraphics[width=4.5cm]{t_x.eps}\includegraphics[width=4.5cm]{t_y.eps} \includegraphics[width=4.5cm]{t_z.eps}
%\end{tabular}
\includegraphics[width=4cm]{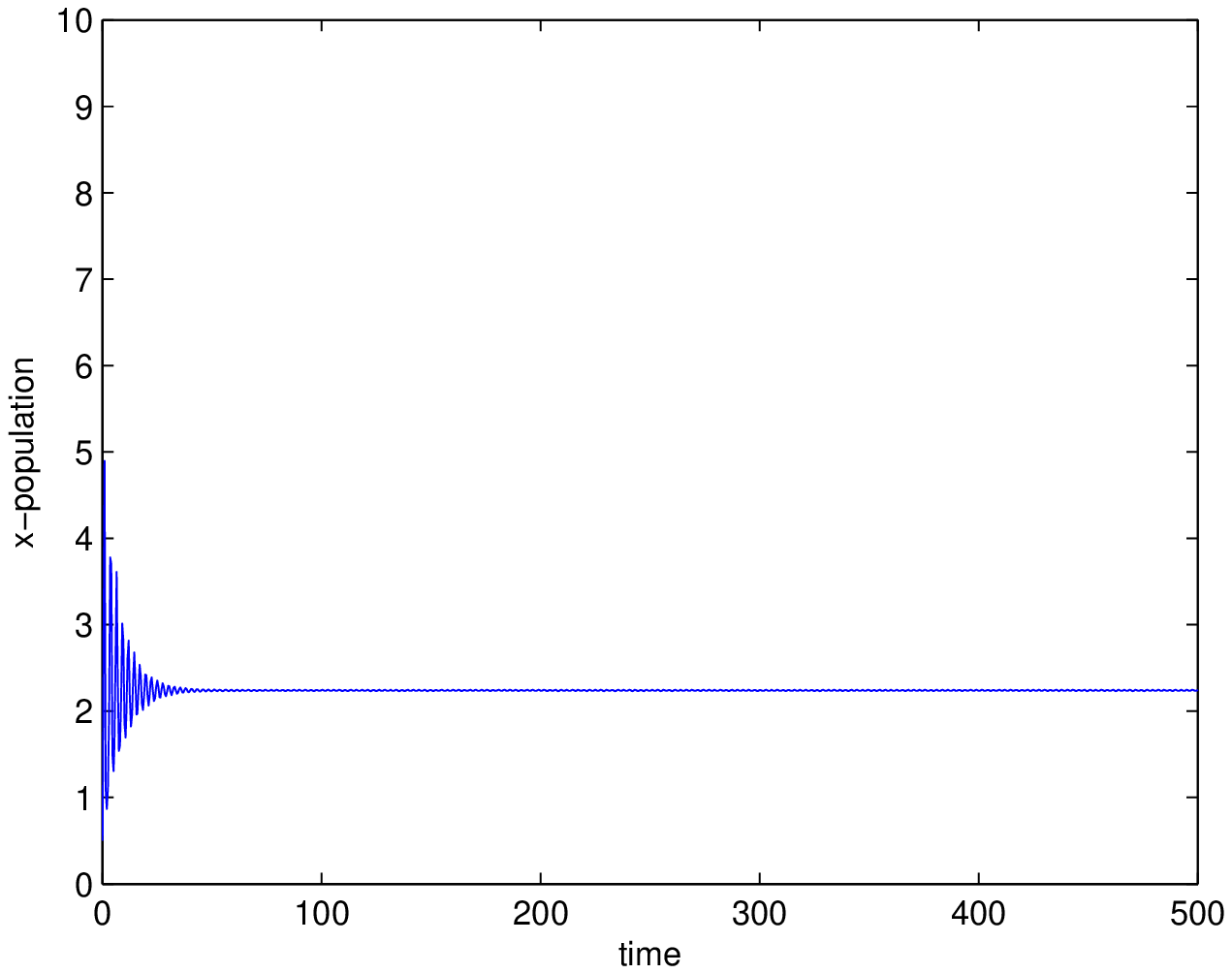}\ \includegraphics[width=4cm]{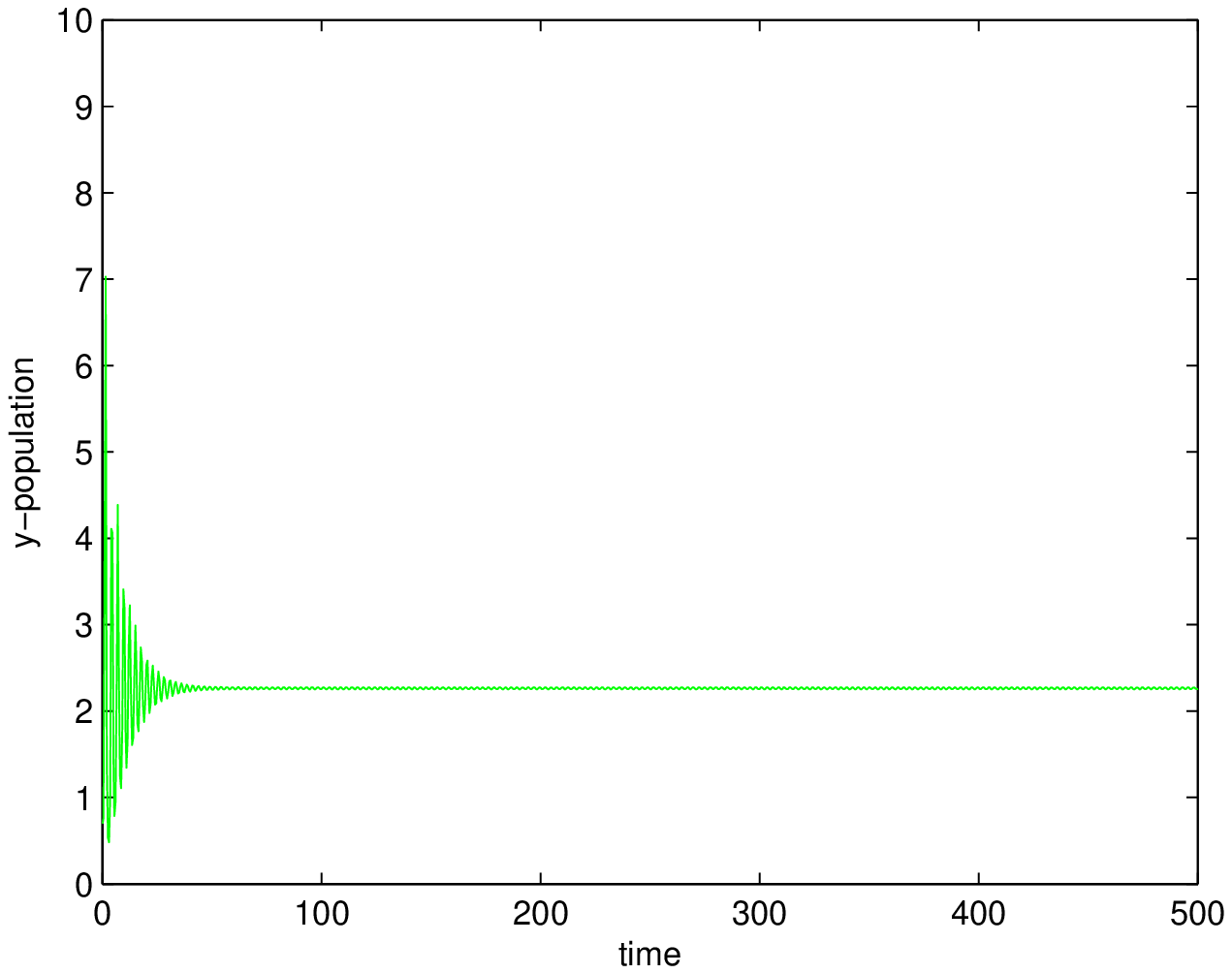}\ \includegraphics[width=4cm]{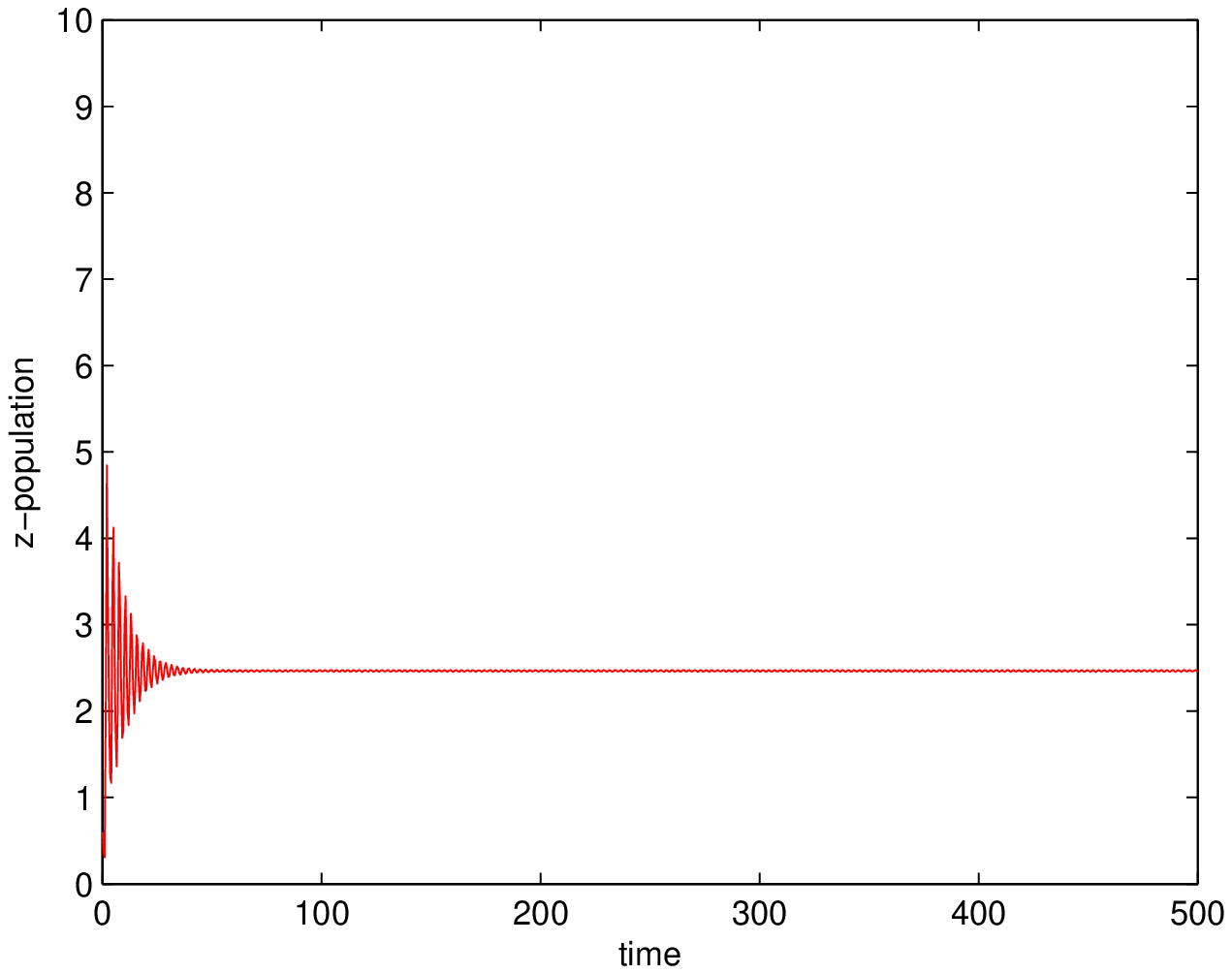}
\caption{This figure shows that the coexistence equilibrium can be stably achieved. Each frame contains, left to right, the populations
$x$, $y$ and $z$ as functions of time.}
\label{stable_focus}
\end{figure}

The study of the coexistence equilibrium appears to be a very difficult task. Even establishing the existence of the equilibrium amounts
to solve a highly nonlinear algebraic system. An attempt along these lines has been made in \cite{HMV10}, see also \cite{HMV09}, but will not be repeated here.
Rather, numerical simulations show the feasibility and stability of this point, Fig. \ref{stable_focus}. This is obtained using
the parameter set, mostly taken from \cite{Tanabe2005},
\begin{eqnarray}\label{param}
b_1=5, \quad b_2=1, \quad b_3=1.25, \quad c=1, \quad 
a_{11}=0.4, \quad a_{12}=1, \\ \nonumber
a_{21}=1, \quad 
a_{23}=1, \quad a_{32}=1, \quad a_{13}=1.5, \quad a_{31}=0.1.
\end{eqnarray}
In fact, stability hinges on the conditions
\begin{equation}\label{stab_E*}
\displaystyle{{\sigma}_1>0}, \quad
\displaystyle{{\sigma}_3>0}, \quad
\displaystyle{{\sigma}_1{\sigma}_2>{\sigma}_3},
\end{equation}
for the details see the Appendix.
The above parameter choice indeed gives
$\sigma_1=2.4372>0$, $\sigma_2=6.2607$, $\sigma_3=12.5946>0$ and
$\sigma_1 \sigma_2-\sigma_3=2.6637>0$, i.e. the Routh-Hurwitz stability criterion is satisfied. This is further checked by the
eigenvalues of (\ref{EQ:eqn1.2}) which have all negative real parts, namely $-2.1969$ and $-0.1201 \pm 2.3913i$.

Chaotic behavior, as discovered in \cite{Tanabe2005} can be analytically be shown to arise also in this more general model. Again we
refer the reader to the Appendix for the detailed calculations. 

\begin{figure}[ht]
\centering
%\begin{tabular}{ccc}
%\includegraphics[width=4.5cm]{bi_a13_x.eps}& \includegraphics[width=4.5cm]{bi_a13_y.eps} & \includegraphics[width=4.5cm]{bi_a13_z.eps}\\
%\end{tabular}
\includegraphics[width=4cm]{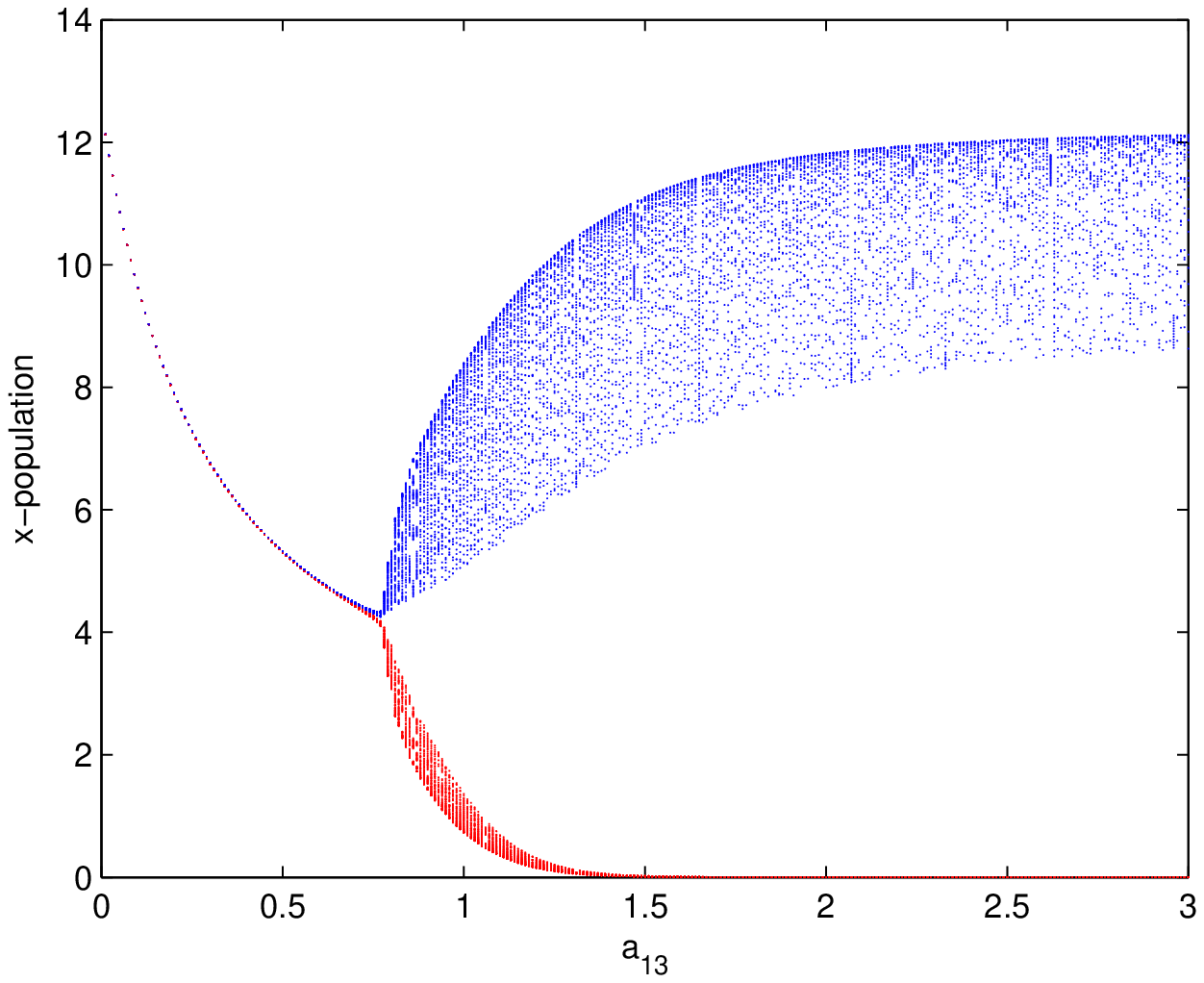}\ \includegraphics[width=4cm]{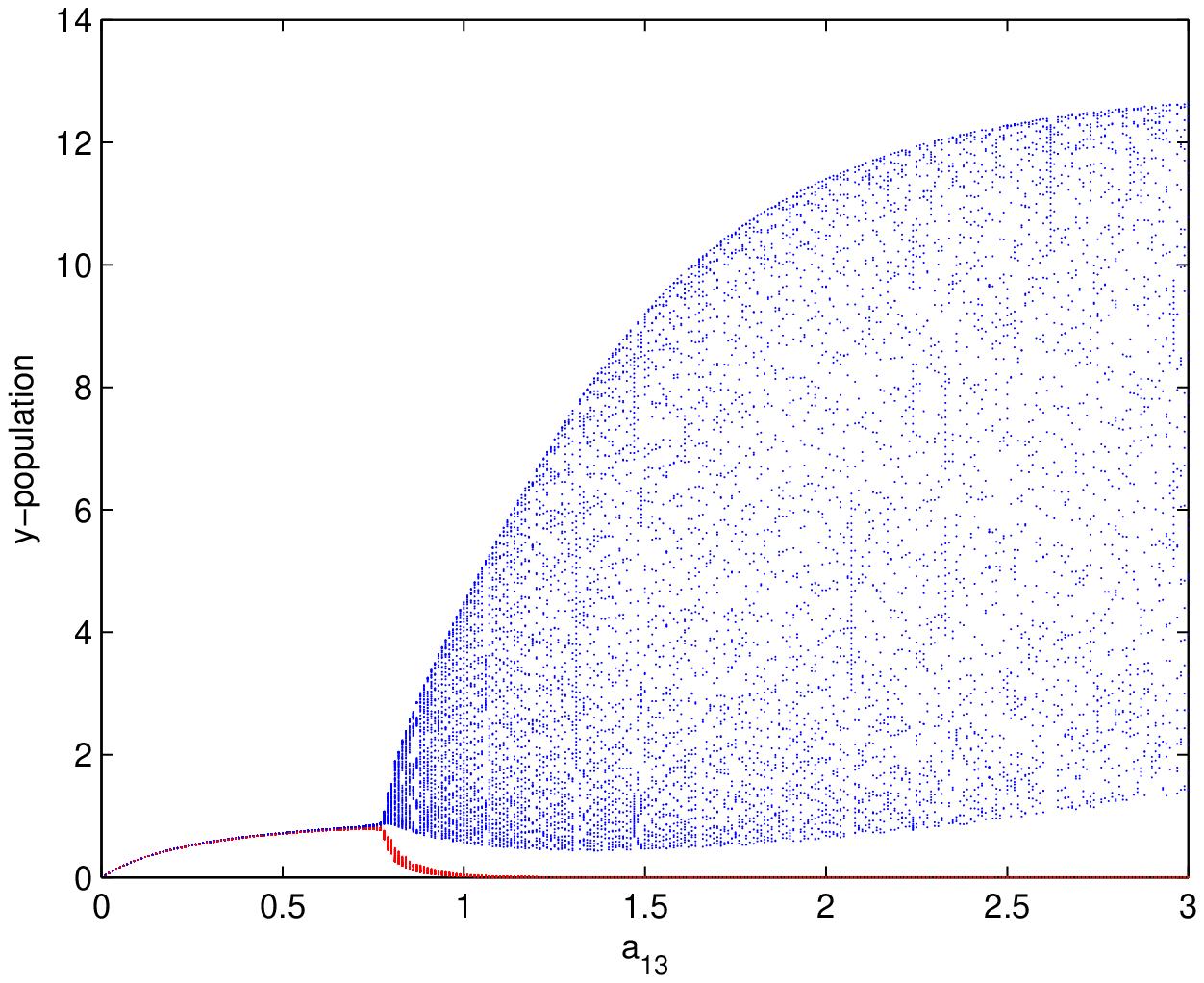}\ \includegraphics[width=4cm]{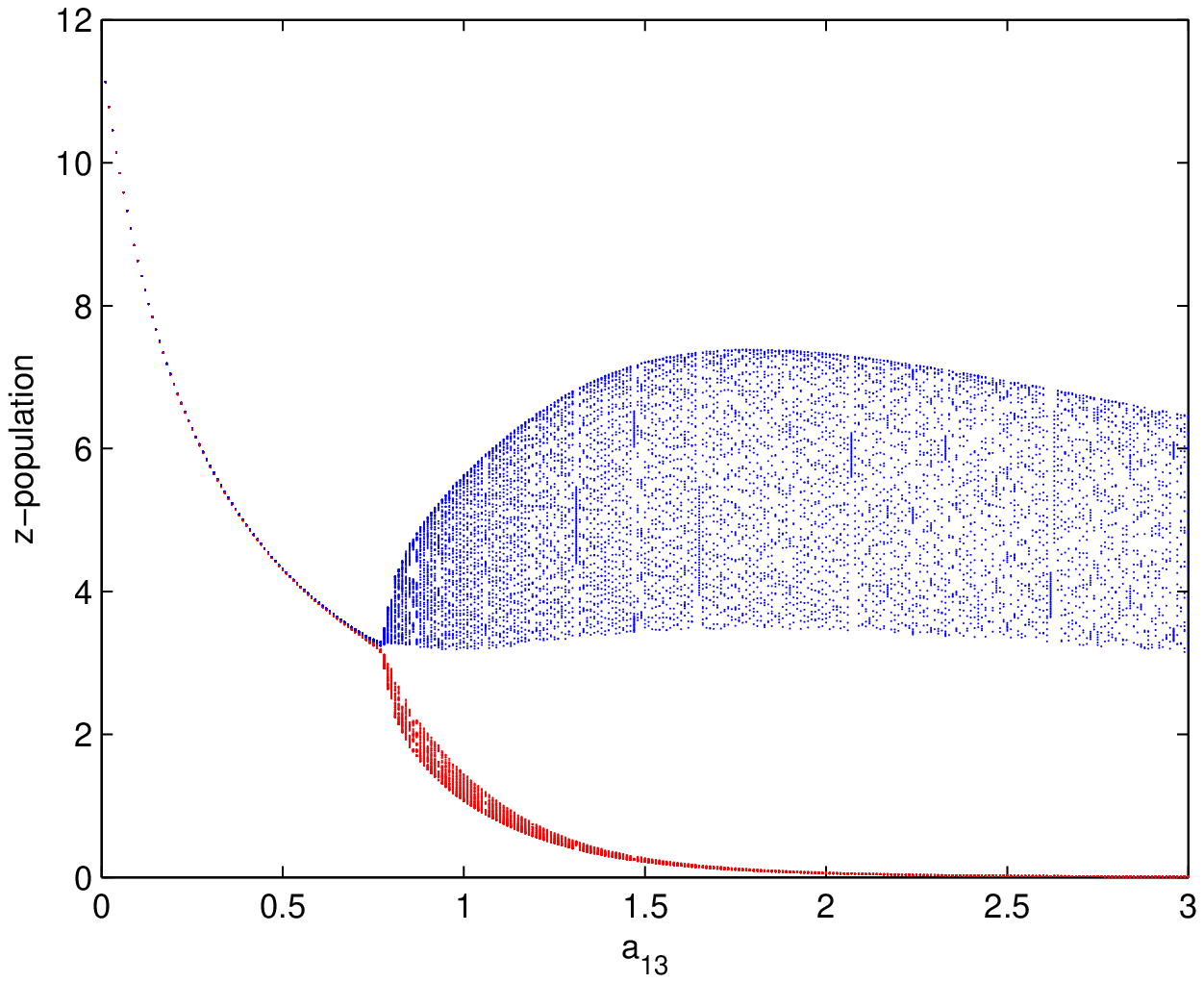}\\
\caption{Bifurcation diagram of the system (\ref{EQ:eqn1.2}) with respect to the bifurcating parameter $a_{13}$.
The other parameter values are the same as those given in (\ref{param}).
Left to right, the diagrams for the populations
$x$, $y$ and $z$ are shown.}
\label{bifurcation_1}
\end{figure}

For the set of parameter values (\ref{param}), the eigenvalues of the characteristic equation of the system
(\ref{EQ:eqn1.2}) are $-1.5765$, $0.2049 \pm 2.4647i$. Therefore we can expect theoretically
the occurrence of Sil'nikov chaos at the coexistence equilibrium ${E}^*$ of the system (\ref{EQ:eqn1.2}).
This theoretical finding is shown
in the bifurcation diagram of Fig. \ref{bifurcation_1}), keeping $a_{13}$ free.

Clearly, if the top predator hunting rate on the bottom prey gradually increases,
the system loses stability and becomes chaotic.

\section{Control of chaos}

Our main purpose here is to investigate the role that predators feeding switching possibly has on the system's behavior.
Recalling that (\ref{EQ:eqn1.2}) exhibits a continuum of models in terms of the switching intensity, we present a simulation
for the bifurcation behavior of the coexistence equilibrium.

\begin{figure}[ht]
\centering
%\begin{tabular}{ccc}
%(a)\includegraphics[width=4.5cm]{bi_c_x.eps} & (b)\includegraphics[width=4.5cm]{bi_c_y.eps}
%(c)\includegraphics[width=4.5cm]{bi_c_z.eps}\\
%\end{tabular}
\includegraphics[width=4cm]{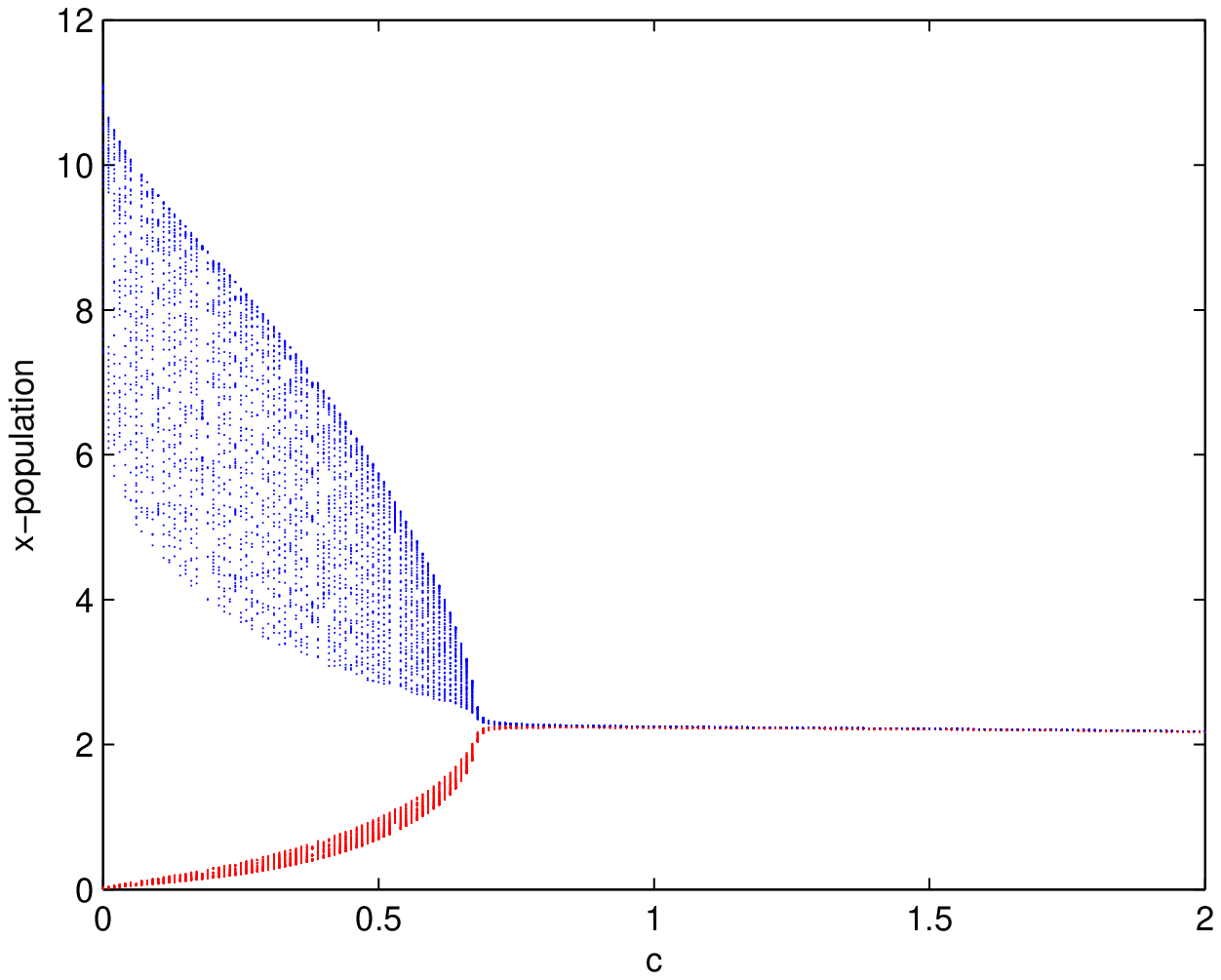} \ \includegraphics[width=4cm]{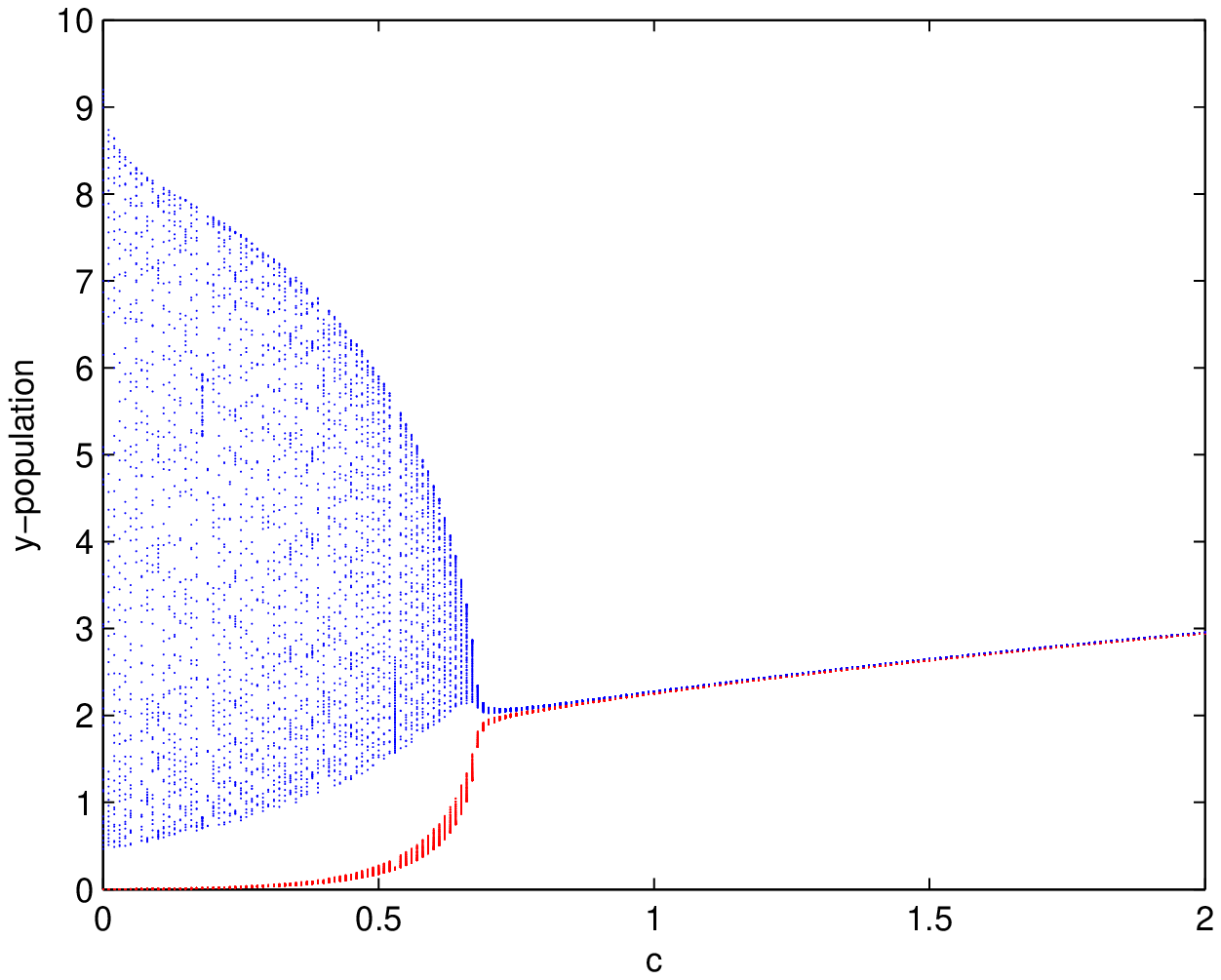}
\ \includegraphics[width=4cm]{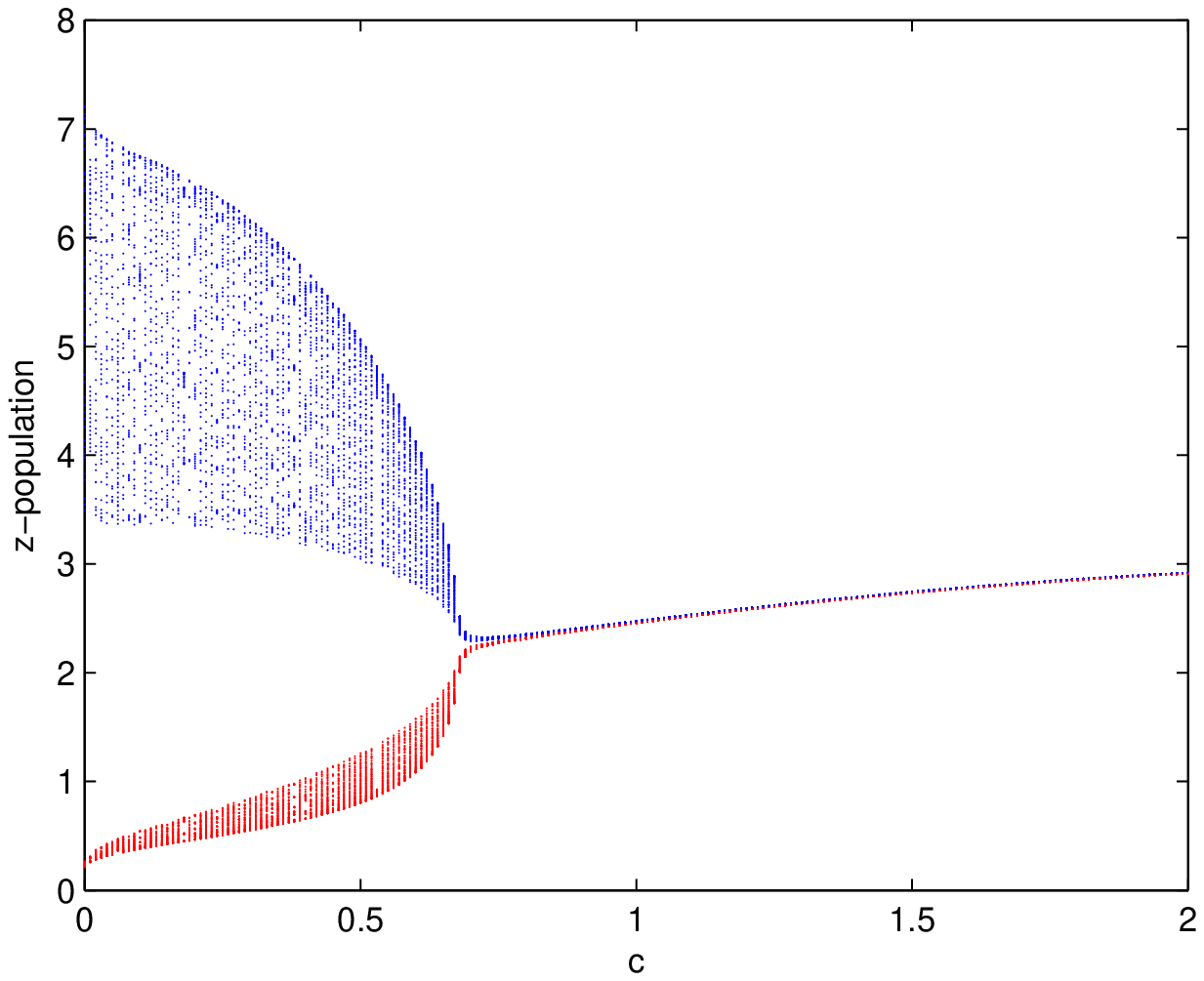}
\caption{Bifurcation diagram of the system (\ref{EQ:eqn1.2}) in terms of the bifurcation parameter $c$.
Here $a_{13}=1.5$ and all the other parameter values are those given in (\ref{param}).}
\label{bifurcation_2}
\end{figure}

The bifurcation diagram of Fig. \ref{bifurcation_2} shows that the system (\ref{EQ:eqn1.2}) switches its
behavior from chaos to a stable focus if the switching intensity increases above a threshold value, which
is found to be approximately $c\approx 0.7$. For $0<c<0.7 $ the model (\ref{EQ:eqn1.2}) exhibits chaotic oscillations
while past that value the coexistence equilibrium is stably achieved.

Our analytical and numerical results indicate thus that if the switching intensity is above a threshold value, then the omnivory system
from the chaotic oscillations achieves a stable behavior, thereby suggesting that predator
switching can in fact enhance the persistence and stability of an omnivory system.
A possible explanation could be provided as follows.
When a particular prey species population declines,
possibly partly owing to the predation, the predator instead of wasting time in looking for this scarce resource
switches its attacks to the next more abundant available prey.
In this way, none of the prey populations is drastically reduced, because hunting stops or reduces drastically when its numbers decline,
nor any prey population is allowed to
become too much abundant, because in this case the predation will resume.
Such mechanism mediates the prey predator interaction in an omnivory system and thus enhances
the stability and the persistence of the system itself.

\section*{Appendix}

\subsection*{Equilbria}
The instability of the origin can be shown by following the method used in \cite{Arino2004}, wo which we refer the reader without
going into details any further here.

The stability conditions stated in the main body of the paper follow from the
calculation of the eigenvalues of the Jacobian matrix evaluated at each equilibria.

In particular for the coexistence equilibrium we find that conditions (\ref{stab_E*}) must hold, where
\begin{eqnarray*}
{\sigma}_1=V_1+V_5,\quad {\sigma}_2=V_1V_5-V_2V_4+V_3V_7+V_6V_8, \\
{\sigma}_3=V_1V_6V_8+V_2V_6V_7+V_3V_4V_8+V_3V_5V_7
\end{eqnarray*}
with
\begin{eqnarray*}
\displaystyle{V_1=a_{11}x^*+\frac{a_{13}cx^*y^*z^*}{(x^*+cy^*)^2} }, \quad
\displaystyle{V_2=-a_{12}x^*+\frac{a_{13}c{x^*}^2z^*}{(x^*+cy^*)^2} },
\displaystyle{V_3=\frac{a_{13}{x^*}^2}{(x^*+cy^*)} }, \\
\displaystyle{V_4=a_{21}y+\frac{a_{23}c{y^*}^2z^*}{(y^*+cx^*)^2} }, \quad
\displaystyle{V_5=\frac{a_{23}cx^*y^*z^*}{(y^*+cx^*)^2} }, \quad
\displaystyle{V_6=\frac{a_{23}{y^*}^2}{(y^*+cx^*)} }, \\
\displaystyle{V_7=\frac{a_{31}{x^*}^2z^*+2a_{31}cx^*y^*z^*}{(x^*+cy^*)^2}-\frac{a_{32}c{y^*}^2z^*}{(y^*+cx^*)^2}}, \\
\displaystyle{V_8=-\frac{a_{31}c{x^*}^2z^*}{(x^*+cy^*)^2}+\frac{a_{32}{y^*}^2z^*+2a_{32}cx^*y^*z^*}{(y^*+cx^*)^2}}.
\end{eqnarray*}

\subsection*{Chaotic behavior}

The characteristic equation of the Jacobian matrix of the system (\ref{EQ:eqn1.2}) at the coexistence equilibrium
$E^*$  is  $\lambda^3+3A_1\lambda^2+3A_2\lambda+A_3=0$,
where
$$
3A_1={\sigma}_1, \quad 3A_2={\sigma}_2, \quad A_3={\sigma}_3.
$$
Now setting $\lambda=\rho-A_1$ the characteristic equation can be written in the form  $\rho^3+3H\rho+G=0$,
where $H=A_2-{A_1}^2$ and $G=A_3-3A_1A_2+2{A_1}^3$.
The roots of the characteristic equation are then
$$
\lambda_1=R-\frac{H}{R}-A_1, \quad
\lambda_2=R\omega-\frac{H}{R}\omega^2-A_1,\quad
\lambda_3=R\omega^2-\frac{H}{R}\omega-A_1,
$$
where
$$
R=\left(-\frac{G}{2} +\frac{\sqrt{\Delta}}{2}\right)^{\frac{1}{3}}, \quad \Delta=G^2+4H^3, \quad
\omega=-\frac{1}{2}+i\frac{\sqrt{3}}{2}, \quad \omega^2=-\frac{1}{2}-i\frac{\sqrt{3}}{2}.
$$

A system has Sil'nikov chaos if the equilibrium point of the system is a saddle focus and the eigenvalues
$\gamma$, $\alpha\pm i\beta$ satisfy the following conditions \cite{Mandal2010},
$$
i) \ \beta\neq0,\quad ii) \ \gamma\alpha<0,\quad iii) \ |\gamma|>|\alpha|\geq0.
$$
For the system (\ref{EQ:eqn1.2}) the above conditions yield
\begin{eqnarray*}
\Delta>0,\quad  R+\frac{H}{R}\neq0,\quad \left| R-\frac{H}{R}-A1\right|- \frac 12 \left| R-\frac{H}{R}+2A1\right|>0, \\
R-\frac{H}{R}+2A_1>0, {\textrm{ for }} A_3<0, {\textrm{ or }} R-\frac{H}{R}+2A_1<0,
{\textrm{ for }} A_3>0.
\end{eqnarray*}

In the absence of the switching mechanism, i.e. for $c=0$ and for the remaining parameters as given in (\ref{param}), we find
$$
\Delta=76.9751>0, \quad R+\frac{H}{R}=2.8460\neq0, \quad R-\frac{H}{R}+2A_1=-0.4099<0,
$$
as well as
$$
A_3=9.6432>0, \quad |R-\frac{H}{R}-A1|-(1/2)|R-\frac{H}{R}+2A1|=1.3716>0,
$$
so that the Sil'nikov conditions for chaotic behavior hold.

\end{document}